\documentclass[12pt, leqno]{article}
\usepackage{amssymb,amsmath,latexsym,amscd}
\evensidemargin =\oddsidemargin

\font\teneufm=eufm10
\font\seveneufm=eufm7
\font\fiveeufm=eufm5
\newfam\eufmfam
\textfont\eufmfam=\teneufm
\scriptfont\eufmfam=\seveneufm
\scriptscriptfont\eufmfam=\fiveeufm

\let\goth\mathfrak

\def\gg{\goth g}

\def\gz{\goth z}

\def\L{\mathcal L}
\def\GL{\text{\rm GL}}
\def\bAut{\textbf{Aut}}

\def\beq{\begin{equation}}
\def\eeq{\end{equation}}
\def\bea{\begin{eqnarray}}
\def\eea{\end{eqnarray}}
\def\beas{\begin{eqnarray*}}
\def\eeas{\end{eqnarray*}}

\def\cplus{\hbox{$\supset${\raise1.05pt\hbox{\kern -0.55em
${\scriptscriptstyle +}$}}\ }}

\DeclareMathOperator{\Hom}{Hom}
\DeclareMathOperator{\Aut}{Aut}

\DeclareMathOperator{\End}{End}
\DeclareMathOperator{\Id}{Id}

\newtheorem{lemma}[equation]{Lemma}

\newtheorem{proposition}[equation]{Proposition}

\newtheorem{example}[equation]{Example}

\newtheorem{remark}[equation]{Remark}

\def\bT{\rm \bf  T}

\title{Descent constructions for central extensions of infinite dimensional Lie
algebras}
\author{Arturo Pianzola\thanks{Supported by the NSERC Discovery Grant Program. The author also wishes to thank
the Instituto Argentino de Matem\'atica for their hospitality.},
Daniel Prelat\thanks{Supported by a Research Grant from Universidad
CAECE.} and Jie Sun \,}
\date{}
\begin{document}
\maketitle

\begin{abstract}
We use Galois descent to construct central extensions of twisted
forms of split simple Lie algebras over rings. These types of
algebras arise naturally in the construction of Extended Affine
Lie Algebras. The construction also gives information about the
structure of the group of automorphisms of such algebras.
\end{abstract}


\bigskip

\section{Introduction}

Throughout $k$ will denote a field of characteristic $0$, and $\gg$
a finite dimensional split simple Lie algebra over $k$.

Given an associative, unital, commutative $k$-algebra $R$, we
consider the (in general infinite dimensional) $k$-Lie algebra
$\gg_R = \gg \otimes_k R$. The case of $R = k[t^{\pm 1}]$ arises in
the untwisted affine Kac-Moody theory, whereof one knows that the
``correct" object to study from the representation point of view is
not $\gg_R$ itself, but rather its universal central extension.

In the Kac-Moody case, the universal central extension is
one-dimensional. This is not so for the ``higher nullity" toroidal
algebras corresponding to $R = k[t_1^{\pm 1},...,t_n^{\pm 1}]$ with
$n > 1$. The universal central extensions are in these cases
infinite dimensional, and one knows of  many interesting  central
extensions of $\gg_R$ which are not universal (see \cite{MRY} and
\cite{EF}).

In this short note, we look at Lie algebras $\cal L$ which are
twisted forms of $\gg_R$. These algebras appear naturally in the
study of Extended Affine Lie Algebras, and present a beautiful
bridge between Infinite Dimensional Lie Theory and Galois Cohomology
(\cite{AABGP}, \cite{ABFP}, \cite{P2}, and \cite{GP1}).  The purpose
of this short note is to give natural constructions for central
extensions of such algebras by descent methods, and to study their
group of automorphisms.

\section{Some generalities on central extensions}
Let  $\L$ be a Lie algebra over $k$, and $V$ a $k$-space. Any
cocycle $P \in Z^2(\L,V)$, where $V$ is viewed as a trivial
$\L$-module, leads to a central extension

$$0\longrightarrow V\longrightarrow \L_P \buildrel \pi
\over \longrightarrow \L\longrightarrow 0$$
of $\L$ by $V$ as
follows: As a space $ \L_P = \L\oplus V$, and the bracket $[\, ,
\,]_P$ on $\L_P$ is given by

$$[x+u,y+v]_{P}=[x,y]+P(x,y) \,\, \text{for} \,\, x,y\in \L \,\, \text{and} \,\, u,v\in
V.$$ The isomorphism class of this extension depends only on the
class of $P$ in $H^2(\L,V)$, and this gives in fact a
parametrization of all isomorphism classes of central extensions of
$\L$ by $V$ (see for example \cite{MP} or \cite{We} for details). In
this situation,  we will henceforth naturally identify $\L$ and $V$
with  subspaces of $\L_P$.

\medskip

An automorphism $\theta \in \Aut_{k}(\L)$ is said to {\it lift to}
$\L_P$, if there exists an element $\theta_P\in \Aut_{k}(\L_P)$ for
which the following diagram commutes.

$$
\begin{CD}
 \L_{P} @>{\pi}>> \L  \\
 @V{\theta_{P}}VV @V{\theta}VV\\
\L_{P} @>{\pi}>> \L \\
\end{CD}
$$
 We then say that $\theta_P$ is a {\it lift} of $\theta$.

\begin{remark} {\rm By definition, $\theta_P$ stabilizes
$V$; thereof inducing an element of $\GL_k(V)$. By contrast,
$\theta_P$ {\it need not} stabilize the subspace $\L$ of $\L_P$.

 Note that $\theta_P(x) - \theta(x) \in V$ for all $x \in \L$. Since
 $V$ lies inside the centre $\gz(\L_P)$ of $\L_P$, we get the useful
 equality}
\begin{equation}
\label{useful} \theta_P([x,y]_P) = [\theta(x),\theta(y)]_P \,\,
\text{\rm for all} \, \, x,y \in
 \L.
\end{equation}
\end{remark}
\begin{lemma}\label{lift}
Let $\delta : \Hom_k(\L,V) \rightarrow Z^2(\L,V)$ be the coboundary
map, i.e. $\delta(\gamma)(x,y) = -\gamma([x,y])$. For $\theta \in
\Aut_{k}(\L)$ and $P \in Z^2(\L,V)$, the following conditions are
equivalent.

(1) $\theta$  lifts to $\L_P$.

(2) There exists $\gamma \in \Hom_k(\L,V)$ and $\mu \in \GL_k(V)$
such that $\mu \circ P-P\circ(\theta\times\theta)=\delta(\gamma)$.

In particular, for a lift $\theta_P$ of $\theta$ to exist, it is
necessary and sufficient that there exists $\mu \in \GL_k(V)$ for
which, under the natural right action of $\GL_k(V) \times
\Aut_k(\L)$ on $H^2(\L,V)$, the element $(\mu, \theta)$ fixes the
class \text{\rm [}$P$\text{\rm ]}. If this is the case, the lift
$\theta_P$ can be chosen so that its restriction to $V$ coincides
with $\mu$.

\end{lemma}
{\it Proof.} {(1)$\Rightarrow$(2)} Let us denote the restriction
 of $\theta_P$ to $V$ by $\mu$. Define
$\gamma:\L\rightarrow V$ by $\gamma : x \mapsto
\theta_P(x)-\theta(x)$. Then
$\theta_P(x+v)=\theta(x)+\gamma(x)+\mu(v)$ for all $x \in \L$ and
$v \in V$. For all $x$ and $y$ in $\L,$ we have
\begin{eqnarray*}
& &(\mu \circ P-P\circ(\theta\times\theta))(x,y)=\mu(
P(x,y))-P(\theta(x),\theta(y))\\
&=&\mu([x,y]_P-[x,y])-([\theta(x),\theta(y)]_P -[\theta(x),\theta(y)])\\
&=&\theta_P([x,y]_P-[x,y])-[\theta(x),\theta(y)]_P+[\theta(x),\theta(y)]\\
&=&\theta_P([x,y]_P)-\theta_P([x,y])-[\theta(x),\theta(y)]_P+[\theta(x),\theta(y)]\\
&=&\theta_P([x,y]_P)-(\theta([x,y])+\gamma([x,y]))-[\theta(x),\theta(y)]_P+[\theta(x),\theta(y)]\ \\
&=&[\theta(x),\theta(y)]_P-\gamma([x,y])-[\theta(x),\theta(y)]_P \ \text {(by \ref{useful} above)}\\
&=&-\gamma([x,y])=\delta(\gamma)(x,y).
\end{eqnarray*}

{(2)$\Rightarrow$(1)} Define $\theta_P \in \End_k(\L_P)$ by
\begin{equation}
\label{definition} \theta_P(x + v)=\theta(x)+\gamma(x)+\mu(v)
\end{equation}
for all $x \in \cal L $ and $v \in V$. Then $\theta_P$ is bijective;
its inverse being given by
\newline $\theta_P^{-1}(x+v)=\theta^{-1}(x)-\mu^{-1}\gamma(\theta^{-1}(x))+\mu^{-1}(v)$.
That $\theta_P$ is a Lie algebra homomorphism is straightforward.
Indeed,
\begin{eqnarray*}
& &\theta_P[x+u,y+v]_P=\theta_P([x,y]+P(x,y))\\
&=&\theta([x,y])+\gamma([x,y])+\mu \circ P(x,y)\\
&=&\theta([x,y])+(\mu \circ P-\delta(\gamma))(x,y)\\
&=&[\theta(x),\theta(y)]+P(\theta(x),\theta(y))\\
&=&[\theta(x),\theta(y)]_P=[\theta_P(x+u),\theta_P(y+v)]_P.
\end{eqnarray*}

Since the action of $\Aut_k(\L) \times \GL_k(V)$ on $Z^2(\L,V)$ in
question is given by $P^{(\mu, \theta)} = \mu^{-1} \circ
P\circ(\theta\times\theta)$, the final assertion is
clear.\footnote{This fact is a Lie algebra version of a well known
result in group theory. For a much more general discussion of the
automorphism group of non-abelian extensions, the reader can consult
\cite{Nb} (specially Theorem B2 and its Corollary).} \hfill
$\square$

\bigskip

For future use, we recall the following fundamental fact.
\begin{proposition} \label{fundamental}
Let $\L$ be a perfect Lie algebra over $k$. Then

(1) There exists a (unique up to isomorphism) universal central
extension
$$
0 \longrightarrow V \longrightarrow \widehat{\L} \buildrel \pi \over
\longrightarrow \L \longrightarrow 0.
$$

(2) If  $\L$ is centreless, the centre $\gz(\widehat{\L})$ of
$\widehat{\L}$ is precisely the kernel $V$ of the projection
homomorphism $\pi : \widehat{\L} \rightarrow \L$ above. Furthermore,
the canonical map $Aut_{k}(\widehat{\L})\rightarrow Aut_{k}(\L)$ is
an isomorphism.

\end{proposition}
{\it Proof.} (1) The existence of an initial object in the category
of central extensions of $\L$ is due to Garland \cite{Grl}. (See
also \cite{Ne}, \cite{MP} and \cite{We} for details).

(2) This result goes back to van der Kallen \cite{vdK}. Other proofs
can be found in \cite{Ne} and \cite{P1}. \hfill $\square$

\begin{remark} \label{universalcocycle}{\rm Assume $\L$ is perfect and centreless. We fix once and for all a
 universal central extension
$0 \longrightarrow V \longrightarrow \widehat{\L} \buildrel \pi
\over \longrightarrow \L \longrightarrow 0 $ (henceforth referred to
as {\it the} universal central extension of $\L$). We will find it
useful at times to think of this extension as being given by a
(fixed in our discussion)  ``universal" cocycle\footnote{This
cocycle is of course not
 unique.} $\widehat{P}$, i.e.
$\widehat{\L} = \L_{\widehat{P}} =  \L \oplus V$. The space $V$ is
then the centre of $\widehat{\L}$, and we write $\widehat{\L} =
 \L \oplus \gz(\widehat{\L})$ to emphasize this point.}
\end{remark}

\begin{lemma}\label{equivalent}
Assume $\L$ is centreless and perfect, and let $0\rightarrow
\gz(\widehat{\L}) \rightarrow\widehat{\L}\rightarrow \L\rightarrow
0$ be its universal central extension.  For an automorphism $\theta
\in Aut_{k}(\L)$, the following conditions are equivalent.

(1) The lift $\widehat{\theta}$ of $\theta$ to $\widehat{\L}$ acts
on the centre of $\widehat{\L}$ by scalar multiplication, i.e.
$\widehat{\theta}|_{\gz(\widehat{\L})}=\lambda \Id$ for some
$\lambda\in k^{\times}$.

(2) $\theta$ lifts to every central quotient of $\widehat{\L}$.

(3) $\theta$ lifts uniquely to every central quotient of
$\widehat{\L}$.

(4) $\theta$ lifts to every central extension of $\L$.
\end{lemma}
{\it Proof.} (1)$\Rightarrow$(2) The lift $\widehat{\theta}$ exists
by Proposition \ref{fundamental}.2. Let $\L_P$ be a central quotient
of $\widehat{\L}$. Then $\L_P\simeq \widehat{\L}/J=(\L\oplus
\gz(\widehat{\L}))/J \simeq \L\oplus \gz(\widehat{\L})/J$. Since
$\widehat{\theta}|_{\gz(\widehat{\L})}=\lambda \Id$ for some
$\lambda\in k^{\times}$, we have $\widehat{\theta}(J)\subset J$. So
$\widehat{\theta}$ induces an automorphism of $\L_P$.

(2) $\Rightarrow$ (3) The point is that $0\rightarrow J\rightarrow
\widehat{\L}\rightarrow \L\oplus \gz(\widehat{\L})/J \rightarrow 0$
is a universal central extension of $\L\oplus \gz(\widehat{\L})/J$.
By Proposition \ref{fundamental}.2 then, any two lifts of $\theta$
to $\L\oplus \gz(\widehat{\L})/J$ must coincide (since they both
yield $\widehat{\theta}$ when lifted to $\widehat{\L}$).

 (3)$\Rightarrow$(4) There is no loss of generality in assuming that the central
 extension of $\L$ is given by a cocycle, i.e. that it is of the form  $0\rightarrow
V\rightarrow\L_P\rightarrow \L\rightarrow 0$ for some $P \in Z^2(\L,
V)$. There exists then a unique Lie algebra homomorphisms $\phi :
\widehat{\L} \rightarrow \L_P$ such that the following diagram
commutes.

$$
\begin{CD}
0 @>>> \gz(\widehat{\L}) @>>> \widehat{\L} @>>> \L  @>>> 0\\
&&@V{{\phi_{|}}_{\gz(\widehat{\L})} }VV @V{\phi}VV @V{id_{\L}}VV\\
0 @>>> V @>>> \L_{P} @>>> \L @>>> 0\\
\end{CD}
$$
Let $J = \ker \phi$. Then $J \subset \gz(\widehat{\L})$, and the
central quotient $\L\oplus \gz(\widehat{\L})/J$ corresponds to the
cocycle $Q$  obtained by reducing modulo $J$ the  universal cocycle
$\widehat{P}$ chosen in modeling $\widehat{\L}$ (see Remark
\ref{universalcocycle}). We have $\L_P = \L\oplus V = \phi(\L\oplus
\gz(\widehat{\L})/J)\oplus V^{'}$ for some suitable subspace $V^{'}$
of $V$. Let $\theta_{Q}$ be the unique lift of $\theta$ to $\L\oplus
\gz(\widehat{\L})/J$, which we then transfer, via $\phi$, to an
automorphism $\theta_Q'$ of the subalgebra $\phi(\L\oplus
\gz(\widehat{\L})/J)$ of $\L_P$. Then $\theta_P=\theta_Q' +
\Id_{V^{'}}$ is a lift of $\theta$ to $\L_P$.

 (4)$\Rightarrow$(1) Assume
$\widehat{\theta}|_{\gz(\widehat{\L})}\neq\lambda \Id$ for all
$\lambda \in k^{\times}$. Then, there exists a line $J \subset
\gz(\widehat{\L})$ such that $\widehat{\theta}(J)\nsubseteq J$.
Consider the central quotient $\L\oplus \gz(\widehat{\L})/J$. Let
$\theta_P$ be a lift of $\theta$ to $\L\oplus
\gz(\widehat{\L})/J$. As pointed out before, $0\rightarrow
J\rightarrow \widehat{\L}\rightarrow \L\oplus
\gz(\widehat{\L})/J\rightarrow 0$ is the universal central
extension of $\L\oplus \gz(\widehat{\L})/J$. So $\widehat{\theta}$
is also the unique lift of $\theta_P$. This forces
$\widehat{\theta}(J)\subset J$, contrary to our assumption. \hfill
$\square$

\section{The case of $\gg_R$}

Throughout $\gg$ will denote a finite dimensional split simple Lie
algebra over $k$, and $R$  a commutative, associative, unital
$k$-algebra. We view $\gg_R=\gg\otimes _{k}R$ as a Lie algebra {\it
over} $k$ (in general infinite dimensional) by means of the unique
bracket satisfying
\begin{equation}
[x\otimes a,y\otimes b]=[x,y]\otimes ab
\end{equation}
for all $x,y \in \gg$ and $a,b \in R$. Of course $\gg_R$ is also
 naturally an $R$-Lie algebra (which is free of finite rank). It will be at
all times clear which of the two structures is being considered.

Let $(\Omega_{R/k}, d_R)$ be the  $R$-module of K\"{a}hler
differentials of the $k$-algebra $R$. When no confusion is possible,
we will simply write $(\Omega_R, d)$. Following Kassel \cite{Ka}, we
consider the $k$-subspace $dR$ of $\Omega_R$, and the corresponding
quotient map $^{\overline{\,\,\,}} \, : \, \Omega_R \rightarrow
\Omega_R/dR$. We then have a unique  cocycle $\widehat{P} =
\widehat{P}_R \in Z^2(\gg_R , \Omega_R/dR)$ satisfying
\begin{equation}
\widehat{P}(x\otimes a,y\otimes b)=(x|\, y)\overline{adb},
\end{equation}
where $(\,|\,\,)$ denotes the Killing form of $\gg$.

Let $\widehat{\gg_R}$ be the unique Lie algebra over $k$ with
underlying space $\gg_{R}\oplus \Omega_{R}/dR$, and bracket
satisfying
\begin{equation}
[x\otimes a,y\otimes b]_{\widehat{P}} = [x,y]\otimes ab + (x|\,
y)\overline{adb}.
\end{equation}
As the notation suggests,
$$
0 \longrightarrow \Omega_{R}/dR \longrightarrow \widehat{\gg_R}
\buildrel \pi \over \longrightarrow \gg_R  \longrightarrow 0
$$
 is the universal central extension of
$\gg_R$.\footnote{ There are other different realizations of the
universal central extension (see \cite{Ne}, \cite{MP} and
\cite{We} for details on three other different constructions), but
Kassel's model is perfectly suited for our purposes.}

\begin{proposition}\label{Rlift}
Let $\theta \in \Aut_k(\gg_R)$, and let $\widehat{\theta}$ be the
unique lift of $\theta$ to $\widehat{\gg_R}$ (see Proposition
\ref{fundamental}). If $\theta$ is $R$-linear, then
$\widehat{\theta}$ fixes the centre $\Omega_R/dR$ of
$\widehat{\gg_R}$ pointwise. In particular, every $R$-linear
automorphism of $\gg_R$ lifts to every central extension of $\gg_R$.

\end{proposition}
{\it Proof.} For future use, we begin by observing that $[x\otimes
a,x\otimes b]_{\widehat{P}}=(x|x)\overline{adb}$. Let now $\theta\in
\Aut_k(\gg_R)$ be $R$-linear. Fix $x \in \gg$ such that $(x |\, x)
\neq 0$, and write $\theta(x)=\sum_{i}x_{i}\otimes a_{i}$. Then
\begin{eqnarray*}
0&=&[x\otimes ab,x\otimes
1]_{\widehat{P}}=\widehat{\theta}([x\otimes ab,x\otimes
1]_{\widehat{P}})\\
&=&[\theta(x\otimes ab),\theta(x\otimes
1)]_{\widehat{P}} \ (\text{\rm  by \ref{useful})}\\
&=&[\sum_{i}x_{i}\otimes aba_{i},\sum_{j}x_{j}\otimes a_{j}]_{\widehat{P}}\\
&=&\sum_{i,j}[x_{i},x_{j}]\otimes
aba_{i}a_{j}+\sum_{i,j}(x_{i}|x_{j})\overline{aba_{i}da_{j}}.\\
\end{eqnarray*}
Thus
\begin{equation}\label{eq1}
\sum_{i,j}(x_{i}|x_{j})\overline{aba_{i}da_{j}} = 0 =
\sum_{i,j}[x_{i},x_{j}]\otimes aba_{i}a_{j}.
\end{equation}

Since $\theta$ is $R$-linear, it leaves invariant the Killing form
of the $R$-Lie algebra $\gg_R$. We thus have
\begin{equation}\label{eq2}
(x|x)_{\gg} =( x\otimes 1|x\otimes 1)_{\gg_R} =(\theta(x\otimes
1)|\theta(x\otimes 1))_{\gg_R} = \sum_{i,j}(x_{i}|x_{j})a_{i}a_{j}.
\end{equation}

We are now ready to prove the Proposition. By Lemma
\ref{equivalent}, it will suffice to show that $\widehat{\theta}$
fixes $\Omega_{R}/dR$ pointwise. Now,
\begin{eqnarray*}
\widehat{\theta}((x|x)\overline{adb})&=&\widehat{\theta}([x\otimes
a,x\otimes
b]_{\widehat{P}})=[\theta(x\otimes a),\theta(x\otimes b)]_{\widehat{P}}\ (\text{\rm  by \ref{useful})}\\
&=&[\sum_{i}x_{i}\otimes aa_{i},\sum_{j}x_{j}\otimes ba_{j}]_{\widehat{P}}\\
&=&\sum_{i,j}[x_{i},x_{j}]\otimes
aba_{i}a_{j}+\sum_{i,j}(x_{i}|x_{j})\overline{aa_{i}dba_{j}}
=\sum_{i,j}(x_{i}|x_{j})\overline{aa_{i}dba_{j}} \ \text{\rm (by \ref{eq1})} \\
&=&\sum_{i,j}(x_{i}|x_{j})\overline{aa_{i}bda_{j}}+\sum_{i,j}(x_{i}|x_{j})\overline{aa_{i}a_{j}db}
=\sum_{i,j}(x_{i}|x_{j})\overline{aa_{i}a_{j}db}\ \text{\rm (by \ref{eq1})}\\
&=&\overline{\sum_{i,j}(x_{i}|x_{j})a_{i}a_{j}adb}=(x|x)\overline{adb}\
\text{\rm (by \ref{eq2}). }
\end{eqnarray*}
\hfill $\square$

\begin{remark}\label{looplift}
{\rm Each element $\theta \in \Aut_k(R)$ can  naturally be viewed as
an automorphism of the Lie algebra $\gg_R$ by acting on the
$R$--coordinates, namely $\theta(x\otimes r) = x \otimes \,
^{\theta}r$ for all $x \in \gg$ and $r \in R$. The group $\Aut_k(R)$
 acts naturally as well on the space $\Omega_{R}/dR$, so that
$^{\theta}(\overline{adb}) = \overline{^{\theta}a d^{\theta}b}$. A
straightforward calculation shows that the map $\widehat{\theta} \in
\GL_k(\widehat{\gg_R})$ defined by $\widehat{\theta}(y + z) =
\theta(y) + \, ^{\theta}z$ for all $y \in \gg_R$ and $z \in
\Omega_{R}/dR$, is an automorphism of the Lie algebra
$\widehat{\gg_R}$. Thus $\widehat{\theta}$ is the unique lift of
$\theta$ to $\widehat{\gg_R}$ prescribed by Proposition
\ref{fundamental}.2. Note that $\widehat{\theta}$ stabilizes the
subspace $\gg_R$.

 We now make some general observations about the automorphisms of $\gg_R$ that lift
 to a given central quotient of $\widehat{\gg_R}$. Without loss of
 generality, we assume that the central quotient at hand is of the form
$(\gg_R)_P$ for some $P \in Z^2(\gg_R, V)$. Let $\theta \in
\Aut_k(\gg_R)$. Since $\widehat{\gg_R}$ is the universal central
extension of its central quotients, the lift $\theta_P$, if it
exists, is unique (Lemma \ref{equivalent}).
 We have $\Aut_k(\gg_R) = \Aut_R(\gg_R) \rtimes \Aut_k(R)$ (\cite{ABP2} Lemma 4.4. See also \cite{BN} Corollary 2.28).
 By Proposition \ref{Rlift} all
 elements of
 $\Aut_R(\gg_R)$ do lift, so the problem reduces to understanding
 which
 $\theta \in \Aut_k(R)$  admit a lift $\theta_P$ to  $(\gg_R)_P$.
    Since
$\widehat{\theta}$ stabilizes $\gg_R$, the linear map $\gamma$ of
Lemma \ref{lift} vanishes.\footnote{This holds for every central
extension of $\gg_R$, and not just central quotients of
$\widehat{\gg_R}$. We leave the details of this general case to
the reader.} We conclude that $\theta_P$ exists if and only if
there exists a linear automorphism $\mu \in \GL_k(V)$ such that $P
= P^{(\mu, \theta)}$. }
\end{remark}

\begin{example} {\rm  Let $k =  \Bbb C$ and $R = \Bbb C[t_{1}^{\pm1},t_{2}^{\pm1}]$. Fix
$\zeta \in \Bbb{C}$, and consider the one dimensional central
extension $\L_{P_{\zeta}} =\gg_R \oplus \Bbb C c$, with cocycle
$P_{\zeta}$ given by $P_{\zeta}(x\otimes
t_{1}^{m_{1}}t_{2}^{m_{2}},y\otimes t_{1}^{n_{1}}t_{2}^{n_{2}})=(x
|\,y)(m_{1}+\zeta
m_{2})\delta_{m_{1}+n_{1},0}\delta_{m_{2}+n_{2},0}c$ (see
\cite{EF}). We illustrate how our methods can be used to describe
the group of automorphisms of this algebra.

As explained in the previous Remark, $\Aut_{\Bbb C} (\gg_R) =
\Aut_{R}(\gg_R)\rtimes \Aut_{\Bbb C}(R)$, all $R$-linear
automorphisms of $\gg_R$ lift to $\L_{P_{\zeta}}$, and we are down
to  understanding which elements of $\Aut_{\Bbb C}(R)$ can be lifted
to $\Aut_{\Bbb C}(\L_{P_{\zeta}})$.

Each $\theta\in \Aut_{\Bbb C}(R)$ is given by
$\theta(t_{1})=\lambda_{1}t_{1}^{p_{1}}t_{2}^{p_{2}}$ and
$\theta(t_{2})=\lambda_{2}t_{1}^{q_{1}}t_{2}^{q_{2}}$ for some
$\left(
\begin{array}{cc}
p_{1}&p_{2}\\
q_{1}&q_{2}
\end{array}
\right)\in \GL_{2}(\Bbb{Z})$, and $\lambda_{1},\lambda_{2}\in {\Bbb
C}^{\times}$. The natural copy of the torus $\bT = {\Bbb
C}^{\times}\times {\Bbb C} ^{\times}$ inside $\Aut_{\Bbb C}(R)$
clearly lifts to $\Aut_{\Bbb C}(\L_{P_{\zeta}})$. We are thus left
with describing the group $\GL_{2}(\Bbb{Z})_{\zeta}$ consisting of
elements of  $\GL_{2}(\Bbb{Z})$ that admit a lift to $\Aut_{\Bbb C}
(\L_{P_{\zeta}})$.  Then  $\Aut_{\Bbb C}(\L_{P_{\zeta}})\simeq
\Aut_{R}(\gg_R)\rtimes \big(({\Bbb C}^{\times}\times {\Bbb C}
^{\times})\rtimes \GL_{2}(\Bbb{Z})_{\zeta}\big)$. Note also that
since $\text{\rm Pic}(R) =1$, the structure of the group
$\Aut_{R}(\gg_R)$ is very well understood \cite{P1}.

The $\GL_{2}(\Bbb{Z})_{\zeta}$ form an interesting 1-parameter
family of subgroups of $\GL_{2}(\Bbb{Z})$ that we now describe. Let
$\theta \in \GL_2(\Bbb{Z})$. Then by Lemma \ref{lift}, $\theta \in
\GL_{2}(\Bbb{Z})_{\zeta}$ if and only if there exists
$\gamma:\gg_R\rightarrow {\Bbb C}c$ and $\mu \in \GL_{\Bbb C}({\Bbb
C}c) \simeq {\Bbb C}^{\times}$, such that
\begin{equation}\label{eqexample}
\gamma([x\otimes a,y\otimes b]) = (P\circ(\theta\times\theta)-\mu
P)(x\otimes a,y\otimes b)
\end{equation}
 for all $x, y \in \gg$, and
all $a, b \in R$ (in fact $\gamma = 0$, as explained in Remark
\ref{looplift}) . Choose $x,y\in \gg$ with $(x|\,y)\neq0$. Since
in $\gg_R$ we have $[x\otimes 1,y\otimes 1]=[x\otimes
t_{1},y\otimes t_{1}^{-1}] = [x\otimes t_{2},y\otimes
t_{2}^{-1}]$,  a straightforward computation based on
(\ref{eqexample}) yields

\begin{equation}\label{line}
\theta \left(\begin{array}{r}
1 \\
\zeta
\end{array}\right)
:= \left(\begin{array}{rr}
p_1 & p_2  \\
q_1 & q_2
\end{array}\right)
\left(\begin{array}{r}
1 \\
\zeta
\end{array}\right)
= \mu \left(\begin{array}{rr}
1  \\
\zeta
\end{array}\right).
\end{equation}
The group $\GL_{2}(\Bbb{Z})_{\zeta}$ could thus be trivial, finite,
or even infinite, depending on some arithmetical properties of the
number $\zeta.$ For example if $\zeta^2 = -1$, then
$\GL_{2}(\Bbb{Z})_{\zeta}$ is a cyclic group of order $4$, generated
by the element $\sigma$ for which $\sigma(t_1) = t_2$ and
$\sigma(t_2) = t_1^{-1}$.

 The one-parameter family $\GL_{2}(\Bbb{Z})_{\zeta}$ has the
following interesting geometric interpretation (which was suggested
to us by the referee).  By the universal nature of $\widehat{\gg_R}
= \gg_R \oplus \Omega_R/dR$, we can identify the space $Z^2(\gg_R,
\Bbb C)$ of cocycles with $\Hom(\Omega_R/dR, \Bbb C) =
(\Omega_R/dR)^*$. Furthermore, this identification is compatible
with the respective actions of the group $\Aut_k(R)$.

 The action of the torus $ {\rm \bf T}$ on $\Omega_R/dR$ is diagonalizable, and
the fixed point space $V = (\Omega_R/dR)^{\rm \bf T}$ is two
dimensional with basis $\{ \overline {t_1^{-1}dt_1},
\overline{t_2^{-1}dt_2} \}$. The cocycle $P_{\zeta}$ is
$\bT$--invariant, and corresponds to the linear function $F_{\zeta}
\in (\Omega_R/dR)^*$ which  maps $\overline {t_1^{-1}dt_1} \mapsto
1$, $\overline {t_2^{-1}dt_2} \mapsto \zeta$, and vanishes on all
other weight spaces of $\bT$ on $\Omega_R/dR$. We can thus identify
$F_{\zeta}$ with an element of $V^*$. The action of $\GL_2(\Bbb Z)
\subset \Aut_{\Bbb C}(R)$ on $\Omega_R/dR$ stabilizes $V$, and is
nothing but left multiplication with respect to the chosen basis
above.

Let $\theta \in \GL_2(\Bbb Z)$. Then $\theta$ lifts to
$\L_{P_{\zeta}}$ if and only if $P_{\zeta} = P_{\zeta}^{(\mu,
\theta)}$ for some $\mu \in \Bbb C^{\times}$ (Remark
\ref{looplift}). With the above interpretation, this is equivalent
to $\theta$ stabilizing the line $\Bbb C F_{\zeta} \subset V^*$;
which is precisely equation (\ref{line}) (after one identifies $V$
with $V^*$ via our choice of basis).}

\end{example}
 \section{The case of twisted forms of $\gg_R$}

We now turn our attention to forms of $\gg_R$ for the flat
topology of $R$, i.e. we look at $R$-Lie algebras $\L$ for which
there exists a faithfully flat and finitely presented extension
$S/R$ for which
\begin{equation}\label{twistedform}
\L\otimes_{R}S \simeq \gg_R\otimes_{R}S \simeq \gg \otimes_k S,
\end{equation}
where the above are isomorphisms of $S$-Lie algebras.

Let $\bAut(\gg)$ be the $k$-algebraic group of automorphisms of
$\gg$. The $R$-group $\bAut(\gg)_R$ obtained by base change is
clearly isomorphic to $\bAut(\gg_R)$. It is an affine, smooth, and
finitely presented group scheme over $R$ whose functor of points is
given by
\begin{equation}
\bAut(\gg_R)(S) =\Aut_S(\gg_R \otimes_{R}S) \simeq \Aut_S(\gg
\otimes_k S).
\end{equation}
 By Grothendieck's theory of descent (see \cite{Mln} and
\cite{SGA3}), we have a natural  bijective map
\begin{equation}
\text{\rm Isomorphism classes of forms of} \,\, \gg_R
\longleftrightarrow H^1_{\acute et}\big(R, \bAut(\gg_R)\big).
\end{equation}

In the case when $R = k[{t_1}^{\pm 1},...,t_n^{\pm 1}]$, the class
of algebras on the left plays an important role in modern infinite
dimensional Lie theory. For $n=1$ the forms in question are nothing
but the affine Kac-Moody algebras (derived modulo their centres. See
\cite{P2}). For general $n$, these algebras yield all the centerless
cores of Extended Affine Lie Algebras (EALA) which are finitely
generated over their centroids \cite{ABFP}.

Neher has shown how to ``build" EALAs out of their centerless cores
(in particular, his methods yield all central extensions of such
cores) \cite{Ne}. We now illustrate how to naturally build central
extensions for twisted forms of $\gg_R$ by descent considerations.
In the case when the descent data corresponds to an EALA, the
resulting algebra is the universal central extension of the
corresponding centreless core.

Henceforth $S/R$ will be  {\it finite} Galois with Galois group $G$
(see \cite{KO}). We will assume that our $\L$ is split by such an
extension\footnote{By the Isotriviality Theorem of \cite{GP1}, this
assumption is superflous for $R = k[{t_1}^{\pm 1},...,t_n^{\pm
1}]$.}, i.e.
\begin{equation}
\L\otimes_{R}S \simeq \gg\otimes_k S
\end{equation}
as $S$-Lie algebras. The descent data corresponding to $\L$, which a
priori is an element of $\bAut(\gg)(S \otimes_R S),$ can now be
thought as being given by
 a cocycle $u \in Z^{1}(G,\Aut_{S}(\gg_ S))$
(usual non-abelian Galois cohomology), where the group $G$ acts on
$\Aut_{S}(\gg_ S) = \Aut_{S}(\gg\otimes_k S)$ via
$^{g}{\theta}=(1\otimes g)\circ\theta\circ(1\otimes g^{-1})$.

As above, we let $(\Omega_{S},d)$ be the module of K\"{a}hler
differentials of $S/k$. The Galois group $G$ acts naturally both on
$\Omega_S$ and on the quotient $k$-space $\Omega_{S}/dS$, in such
way that $^{g}({\overline{sdt}})=\overline{^{g}sd^{g}t}$. This leads
to an action of $G$ on $\widehat{\gg_{S}}$ for which
$$^{g}{(x\otimes
s+z)}=x\otimes\,^{g}{s}+\,^{g}{z}$$ for all $x\in \gg$, $s\in S$,
$z\in\Omega_{S}/dS$, and $g \in G$. One verifies immediately that
the resulting maps are automorphisms of the $k$-Lie algebras
$\widehat{\gg_{S}}$. Indeed,
\newline
$ ^{g}{[x_{1}\otimes s_{1}+ z_{1},x_{2}\otimes
s_{2}+z_{2}]_{\widehat{\gg_S}}} = ^{g}{([x_{1},x_{2}]\otimes
s_{1}s_{2}+(x_{1}|\,x_{2})\overline{s_{1}ds_{2}})}=
[x_{1},x_{2}]\otimes
\,^{g}{s_{1}}{}^{g}{s_{2}}+(x_{1}|\,x_{2})\overline{^{g}{s_{1}}d{}^{g}{s_{2}}}
=[x_{1}\otimes \,^{g}{s_{1}}+\,^{g}{z_{1}},x_{2}\otimes
\,^{g}{s_{2}}+\,^{g}{z_{2}}]_{\widehat{\gg_S}}. $ Accordingly, we
 henceforth identify $G$ with a subgroup of
$\Aut_{k}(\widehat{\gg_{S}})$, and let then $G$ act on
$\Aut_{k}(\widehat{\gg_{S}})$ by conjugation, i.e.
$^{g}{\theta}=g\theta g^{-1}$.

\begin{proposition}\label{twistedcentral}
Let $u=({u_{g}})_{g \in G}$ be a cocycle in
$Z^{1}\big(G,\Aut_{S}(\gg_S)\big)$. Then

(1) $\widehat{u}=(\widehat{u}_{g})_{g \in G}$ is a cocycle in
$Z^{1}(G,\Aut_{k}(\widehat{\gg_S}))$.

(2) $\L_{\widehat{u}}=\{x\in\widehat{\gg_S}:\
\widehat{u}_{g}{}^{g}{x}=x\ \, \text{\rm for all}\, g\in G\}$ is a
central extension of the descended algebra $\L_u$ corresponding to
$u$.

(3) There exist canonical isomorphisms $\gz(\L_{\widehat{u}}) \simeq
(\Omega_{S}/dS)^G \simeq \Omega_R/dR$.
\end{proposition}
{\it Proof.} (1)
$\widehat{u}_{g_{1}}{}^{g_{1}}\widehat{u}_{g_{2}}$ is clearly a
lift to $\Aut_{k}(\widehat{\gg_S})$ of
$u_{g_{1}}{}^{g_{1}}u_{g_{2}}=u_{g_{1}g_{2}} \in \Aut_k(\gg_S)$.
By uniqueness (Proposition \ref{fundamental}.2), we have
$\widehat{u}_{g_{1}g_{2}}=\widehat{u}_{g_{1}}{}^{g_{1}}\widehat{u}_{g_{2}}$.

(2) It is clear that $\L_{\widehat{u}}$ is a $k$-subalgebra of
$\widehat{\gg_{S}}$. Moreover if $x\in \gg_S$ and
$\omega\in\Omega_{S}$ are such that
$x+\overline{\omega}\in\L_{\widehat{u}}$, then
$$x+\overline{\omega}=\widehat{u}_{g}{}^{g}{(x+\overline{\omega})}=\widehat{u}_{g}(^{g}{x}+^{g}{\overline{\omega}})
=u_{g}{}^{g}{x}+(\widehat{u}_{g}-u_{g})(^{g}{x})+^{g}{\overline{\omega}}$$
(this last equality by Proposition \ref{Rlift} applied to $\gg_S$).
Since $(\widehat{u}_{g}-u_{g})(\gg_S)\subset\Omega_{S}/dS$ we get
$x=u_{g}{}^{g}{x}$, hence that $x\in\L_{u}$. Thus
$\pi(\L_{\widehat{u}})\subset\L_{u}$, where

$$
0 \longrightarrow \Omega_{S}/dS \longrightarrow \widehat{\gg_S}
\buildrel \pi \over \longrightarrow  \gg_S \longrightarrow 0
$$
as above, is the universal central extension of $\gg_S$.  Since the
kernel of
$\pi|_{\L_{\widehat{u}}}:\L_{\widehat{u}}\rightarrow\L_{u}$ is
visibly central, the only delicate point is to show that
$\pi(\L_{\widehat{u}})=\L_{u}$.

Fix $x\in\L_{u}$. We must show  the existence of some
$z\in\Omega_{S}/dS$ for which $x+z\in\L_{\widehat{u}}$. For each
$g\in G$, we have $\widehat{u}_{g}{}^{g}{x}=x+x_{g}$ for some
$x_{g}\in\Omega_{S}/dS$. We claim that
 the map $g\mapsto x_{g}$ is a cocycle in $Z^{1}(G,\Omega_{S}/dS)$.
 Indeed,
$$x+x_{g_{1}g_{2}}=\widehat{u}_{g_{1}g_{2}}{}^{{g_{1}g_{2}}}{x}=\widehat{u}_{g_{1}}{}^{g_{1}}\widehat{u}_{g_{2}}{}^{{g_{1}g_{2}}}{x}
=\widehat{u}_{g_{1}}g_{1}\widehat{u}_{g_{2}}{}^{{g_{2}}}{x}=\widehat{u}_{g_{1}}{}^{g_{1}}(x+x_{g_{2}})=
x+x_{g_{1}}+\,^{{g_{1}}}{x_{g_{2}}},
$$ and the claim follows. Given that $G$ is finite and $k$ is of characteristic $0$, we have $H^{1}(G,\Omega_{S}/dS)=0$.
Thus, there exists $z\in\Omega_{S}/dS$ such that $x_{g}=z-\,^{g}{z}$
for all $g\in G$. Then
$\widehat{u}_{g}{}^{g}{(x+z)}=x+x_{g}+\,^{g}{z}=x+z,$ so that
$x+z\in\L_{\widehat{u}}$ as desired.

(3) Because $\L_{u}$ is centerless (see Remark \ref{finalremark}
below), the centre of $\L_{\widehat{u}}$ lies inside
$\Omega_{S}/dS$, hence inside $(\Omega_S/dS)^G$ by the definition of
$\L_{\widehat{u}}$ together with Proposition \ref{Rlift} applied to
$\gg_S$. Thus, under the canonical identification of $\Omega_S/dS$
with a subspace of $\L_{\widehat{u}}$, we have
$\gz(\L_{\widehat{u}}) = (\Omega_S/dS)^G$.

On the other hand $\Omega_{S}/dS\simeq HC_{1}(S)$, where this last
is the cyclic homology of $S/k$. The group $G$ acts naturally on the
 $HC_{1}(S)$, and the canonical isomorphism
$\Omega_{S}/dS\rightarrow HC_{1}(S)$ is $G$-equivariant. Since $S/R$
is Galois, na\"ive descent holds for cyclic homology (\cite{WG}
proposition 3.2). We thus have
$$(\Omega_{S}/dS)^G \simeq HC_{1}(S)^G \simeq HC_{1}(S^G) \simeq HC_{1}(R) \simeq
\Omega_{R}/dR.\  $$ $\hfill \square$

\begin{proposition}\label{twistedlift}
With the above notation, the following conditions are equivalent.

(1) $\L_{\widehat{u}}=\L_{u}\oplus\Omega_{R}/dR$ and $\L_{u}$ is
stable under the action of the Galois group $G.$

(2) $\widehat{u}_{g}(\L_{u})\subset\L_{u}$ for all $g\in G$.

If these conditions hold, then every $\theta\in Aut_{R}(\L_{u})$
lifts to an automorphism $\widehat{\theta}$ of $\L_{\widehat{u}}$
that fixes the centre of $\L_{\widehat{u}}$ pointwise.

\end{proposition}
{\it Proof.} {(1)$\Rightarrow$(2)} Let $x\in\L_{u}$. By assumption
we have $^{g^{-1}}x\in\L_{u} \subset  \L_{\widehat{u}}$ for all
$g\in G$. Thus
$\widehat{u}_{g}(x)=\widehat{u}_{g}{}^{g}(^{g^{-1}}{x})= \,
^{g^{-1}}x\in\L_{u}.$

{(2)$\Rightarrow$(1)} Let $x\in\L_{u}$. From the assumption
$\widehat{u}_{g}(\L_{u})\subset\L_{u}$ we obtain
$\widehat{u}_{g}(x)=u_{g}(x)\in\L_{u}$ for all $g\in G$. Thus
$$^{g^{-1}}x=u_{g^{-1}g}(^{g^{-1}}x)=u_{g^{-1}}{}^{g^{-1}}u_{g}(^{g^{-1}}x)=u_{g^{-1}}{}^{g^{-1}}(u_{g}(x))=u_{g}(x)\in\L_{u}.$$
This shows that $\L_{u}$ is $G$--stable. Now $^{g}x\in\L_{u}$ yields
$\widehat{u}_{g}(^{g}x)=u_{g}(^{g}x)=x.$ Thus implies
$\L_{u}\subset\L_{\widehat{u}}$. By Proposition \ref{Rlift}
$\widehat{u}_g$ fixes $\Omega_S/dS$ pointwise, so we have
$\Omega_{R}/dR\subset\L_{\widehat{u}}$. Thus
$\L_{u}\oplus\Omega_{R}/dR\subset\L_{\widehat{u}}$. Finally, if
$\widehat{x} \in\L_{\widehat{u}}$ and we write $\widehat{x}=x+z$
with $x\in \L_u$ and $z\in\Omega_{S}/dS$ according to Proposition
\ref{twistedcentral}.2, then $\widehat{x}-x=z\in\L_{\widehat{u}}$.
Thus $^{g}{z}=\widehat{u}_{g}{}^{g}{z}=z$. This shows that
$z\in\Omega_{R}/dR$ (see Proposition \ref{twistedcentral}.3).

As for the final assertion, let $\theta_{S}$ be the unique $S$-Lie
automorphism of $\gg_S$ whose restriction to $\L_{u}$ coincides with
$\theta$. Let $\widehat{\theta_{S}}$ be the lift of $\theta_{S}$ to
$\widehat{\gg_S}$. We claim that $\widehat{\theta_{S}}$ stabilizes
$\L_{\widehat{u}}=\L_{u}\oplus\Omega_{R}/dR$. By Proposition
\ref{Rlift} $\widehat{\theta_{S}}$ fixes $\Omega_{R}/dR$ pointwise.
Let $x\in\L_{u}$. Since $\L_{u}$ is perfect (see Remark
\ref{finalremark} below), we can write
$x=\sum_{i}[x_{i},y_{i}]_{\gg_S}$ for some $x_{i},y_{i}\in\L_{u}$.
Thus $x=\sum_{i}[x_{i},y_{i}]_{\widehat{\gg_S}} - z$ for some
$z\in\Omega_{R}/dR$. Then
$\widehat{\theta_{S}}(x)=\sum_{i}[\theta(x_{i}),\theta(y_{i})]_{\widehat{\gg_S}}-z
\in[\L_{u},\L_{u}]_{\widehat{\gg_S}} +
\Omega_{R}/dR\subset\L_{u}\oplus\Omega_{R}/dR=\L_{\widehat{u}}$.
\hfill $\square$

\begin{remark}\label{finalremark}
{\rm  Let $\L$ be a twisted form of $\gg_R$ in the sense of
(\ref{twistedform}) above. By faithfully flat descent
considerations, $\L$ is centreless. Indeed, the centre $\gz(\L)
\subset \L$ is an $R$-submodule of $\L$. Since $S/R$ is faithfully
flat, the map $\gz(\L)\otimes_R S \to \L \otimes_R S \simeq \gg
\otimes S$ is injective. Clearly the image of $\gz(\L)\otimes_R S$
under this map lies inside the centre of $\gg \otimes S$, which is
trivial (as one easily sees by considering a $k$-basis of $S$, and
using the fact that $\gz(\gg) = 0$). Thus $\gz(\L)\otimes_R S =0$,
and therefore $\gz(\L) = 0$ again by faithfull flatness. Similarly
descent considerations (see \S5.1 and \S5.2  of \cite{GP2} for
details) show that $\L$ is perfect, and that the centroid of $\L$,
both as an $R$ and $k$-Lie algebra, coincides with $R$ (acting
faithfully on $\L$ via the module structure).

Assume now that $\L$ is split by a finite Galois extension $S/R$.
Let $G$ be the Galois group of $S/R$. Then $\L \simeq \L_u$ for some
cocycle $u = {(u_{g})}_{{g \in G}}$ as above.
 The $R$-group $\bAut(\L_u)$ is a twisted form of
$\bAut(\gg_R)$ (\S5.4 of \cite{GP2}) , in particular affine, smooth,
and finitely presented. We have $\Aut_R(\L_u) = \bAut(\L_u)(R)$.
Every automorphism of $\L_u$ as a $k$-Lie algebra induces an
automorphism of its centroid. By identifying now the centroid of
$\L_u$  with $R$ as explained above, we obtain the following useful
exact sequence of groups
\begin{equation}
1 \rightarrow \Aut_R(\L_u) \rightarrow \Aut_k(\L_u) \rightarrow
\Aut_k(R).
\end{equation}
 If moreover the descent data for $\L_u$ falls under the assumption of
Proposition \ref{twistedlift} (which includes the EALA
case\footnote{The crucial point is that for EALAs, the algebra
$\L_u$ may be assumed to be a multiloop algebra (\cite{ABFP},
corollary 8.3.5). By a general fact about the nature of multiloop
algebras as forms (see \cite{P2} for loop algebras, and \cite{GP2}
\S6 in general), the cocycle $u$ is a group homomorphism $u : G
\rightarrow \Aut_k(\gg)$. In particular, $u$ is constant (i.e. it
has trivial Galois action). The multiloop algebra $\L_u$ has then a
basis consisting of eigenvectors of the $u_g$'s, and therefore the
second equivalent condition of Proposition \ref{twistedlift}
holds.}), then one also has a very good understanding of the
automorphism group of $\L_{\widehat{u}}$. This group will
undoubtedly play a role in any future work dealing with conjugacy
questions for Extended Affine Lie Algebras (see \cite{P3} and
\cite{P1} for the toroidal case).

Finally, we observe that since $\L_u$ above is perfect, it admits a
universal central extension $\widehat{\L_u}$. By Proposition
\ref{twistedcentral}, there exists a canonical map $\widehat{\L_u}
\rightarrow \L_{\widehat{u}}$. In the case of descent data arising
from EALAs, this map is an isomorphism (work in progress of Neher,
given that $\L_{\widehat{u}} = \L_u \oplus \Omega_R/dR$ by
\cite{ABFP} {\it loc .cit.} as explained above). What happens in
general however, remains an open problem.}
\end{remark}

\noindent{\bf Acknowledgement}. We would like to thank the
 referee for his/her useful comments. We are also grateful to E. Neher for explaining to us his
 construction of central extensions of EALA cores (\cite{N1}, \cite{N2}).

\bigskip
\bigskip

Arturo Pianzola and Jie Sun, Department of Mathematical and
Statistical Sciences, University of Alberta, Edmonton, Alberta T6G
2G1 CANADA.
\medskip

Daniel Prelat, Departamento de Matem\'atica,  Universidad CAECE,
Avenida de Mayo 866, Buenos Aires ARGENTINA.

\bigskip

\end{document}